\documentclass[12pt,a4paper,reqno]{amsart}

\usepackage{anysize}
\usepackage[utf8]{inputenc}
\usepackage{amsmath}
\usepackage{amssymb}
\usepackage{array}

\marginsize{2.5cm}{2.5cm}{3cm}{3cm}

\newtheorem{theorem}{Theorem}[section]

\newtheorem{example}[theorem]{Example}

\newtheorem{remark}[theorem]{Remark}

\begin{document}

\title[Bilateral Ramanujan-like series for $1/\pi^k$ and their congruences]{Bilateral Ramanujan-like series for $1/\pi^k$ \\ and their congruences}
\author{Jesús Guillera}
\address{Department of Mathematics, University of Zaragoza, 50009 Zaragoza, SPAIN}
\email{jguillera@gmail.com}
\date{}

\maketitle

\begin{abstract}
We prove a kind of bilateral semi-terminating series related to Ramanujan-like series for negative powers of $\pi$, and conjecture a type of supercongruences associated to them. We support this conjecture by checking all the cases for many primes. In addition we are able to prove a few of them from some terminating hypergeometric identities. Finally, we make an intriguing observation.
\end{abstract}

\section{Introduction}
A Ramanujan-like series for negative powers of $\pi$ is an identity of the following form:
\begin{equation}\label{rama-series-alg}
\sum_{n=0}^{\infty} \left( \prod_{i=0}^{2m} \frac{(s_i)_n}{(1)_n} \right) \sum_{k=0}^m a_k n^k z_0^n \, = \frac{1}{\pi^m},
\end{equation}
where $2m+1$ is the rank of the series, $z_0$ and  $a_0, a_1...,a_m$ are algebraic numbers, and the $0<s_i<1$ are rational fractions with the following property: If we have a fraction, we also have all irreducible proper fractions with the same denominator. So, it is sure that there is an odd number of fractions $s_i$ that are equal to $1/2$. 
\par These series have a long story that begun in $1914$. In that year, S. Ramanujan discovered the family of series of rank $3$, and found $17$ examples (indicating very briefly the idea of the proof). One of his most famous formulas for $1/\pi$ is
\[
\sum_{n=0}^{\infty} \frac{(\frac12)_n(\frac14)_n(\frac34)_n}{(1)_n^3} (26390n+1103) \left(\frac{1}{99^4} \right)^n = \frac{9801 \sqrt{2}}{4\pi}.
\]
However, they looked enigmatic until J. Borwein and P. Borwein gave the first rigorous proofs of them. Nowadays we know that there are $44$ rational Ramanujan-type series for $1/\pi$, $8$ of which are ``divergent" (Tables \ref{3F2-1} and \ref{3F2-2}). See the excellent surveys \cite{BaBeCh} and \cite{Zudilin-rama-wind}.
\par In $2002$ and $2003$ I discovered the family of rank $5$ finding $7$ examples and proving $3$ of them by the WZ method \cite{PWZ}, one of which is
\[
\sum_{n=0}^{\infty} \frac{(\frac12)_n^5}{(1)_n^5}(820n^2+180n+13) \left( -\frac{1}{1024}\right)^n = \frac{128}{\pi^2}.
\]
In Table \ref{table-5F4} we show all the rational series in the family known up to date (most of them remain unproved). See my joint paper with G. Almkvist \cite{Alm-Gui-String}, and the excellent survey \cite{Zudilin-arith-hyper}.
\par A few formulas (but unproved) of higher rank are also known. Formulas $\#1$ and $\#2$ of Table \ref{table-7F6} were discovered by B. Gourevitch and the author respectively, and the formula $\#3$ of Table \ref{table-7F6} was discovered by Y. Zhao \cite{Zhao-1}. More precisely he gives a reverse convergent version of it for $\pi^4$. Formula $\#1$ of Table \ref{table-9F8}, namely
\[
\sum_{n=0}^{\infty} \frac{(\frac12)_n^7(\frac14)_n(\frac34)_n}{(1)_n^3} (43680n^4+20632n^3+4340n^2+466n+21) \left(\frac{1}{4096} \right)^n = \frac{2048}{\pi^4},
\]
was discovered by J. Cullen, and the formulas $\#2$ and $\#3$ by Y. Zhao \cite{Zhao-1}. More precisely, instead of $\#3$ he gives its reverse, a nice convergent formula for $\zeta(5)$. 
\par In this paper, we only consider rational Ramanujan-like series. They can be written as
\begin{equation}\label{rama-series}
\sum_{n=0}^{\infty} \left( \prod_{i=0}^{2m} \frac{(s_i)_n}{(1)_n} \right) \sum_{k=0}^m a_k n^k z_0^n \, = \frac{v_0\sqrt{(-1)^m \chi_0}}{\pi^m},
\end{equation}
where $z_0$ is a rational, $a_0, a_1,...,a_m$ are positive integers, $v_0$ is an integer, and $\chi_0$ the discriminant of a certain quadratic field (imaginary or real), which is also an integer. Then, we construct and study bilateral series related to (\ref{rama-series}), investigate semi-terminating cases of them, and conjecture associated semicongruences module $p^{2m+1}$. 
\begin{remark} \rm \label{obser-1}
If $|z_0|>1$ then the series on the left side of (\ref{rama-series}) is divergent, but we understand it as the analytic continuation of
\[
\sum_{n=0}^{\infty} \left( \prod_{i=0}^{2m} \frac{(s_i)_n}{(1)_n} \right) \sum_{k=0}^m a_k n^k z^n,
\]
at $z=z_0$.
\end{remark} 
\section{Bilateral series of Ramanujan-type}
Here, we construct bilateral series related to Ramanujan-like series for negative powers of $\pi$. Firs, we define the following functions:
\begin{equation}\label{Anx}
A(n,x)= \left( \prod_{i=0}^{2m} \frac{(s_i)_{n+x}}{(1)_{n+x}} \right) \sum_{k=0}^m a_k (n+x)^k \, z_0^{n+x},
\end{equation}
and
\begin{equation}\label{Bnx}
B(n,x)= z_0^{-x} \left( \prod_{i=0}^{2m} \frac{(1)_x}{(s_i)_x} \right) A(n,x) = \left( \prod_{i=0}^{2m} \frac{(s_i+x)_n}{(1+x)_n} \right) \sum_{k=0}^m a_k (n+x)^k \, z_0^n.
\end{equation}
Then, the following identities for bilateral series hold:
\begin{multline}\label{bilateral-rama}
e^{-i \pi x} \, (\cos \pi x)^{2j+1} \prod_{s_k \neq \frac12} \frac{\cos \pi x - \cos \pi s_k}{1-\cos \pi s_k}  \sum_{n=-\infty}^{\infty} A(n, x) \\ = e^{-i \pi x} \, z_0^x \left( \frac{(\cos \pi x) \left(\frac12\right)_x}{(1)_x} \right)^{2j+1} \prod_{s_k \neq \frac12} \frac{(s_k)_x}{(1)_x} \, \frac{\cos \pi x - \cos \pi s_k}{1-\cos \pi s_k}  \sum_{n=-\infty}^{\infty} B(n, x)
\\ = \frac{v_0\sqrt{(-1)^m \chi_0}}{\pi^m}  \left( 1+ \sum_{k=1}^{m} \left( \alpha_k (\cos 2 \pi k x -1) + \beta_k \sin 2 \pi k x \right) \right),
\end{multline}
where the $\alpha_k$ and $\beta_k$ are constants, that we conjecture are rational. In the case $m=1$ we can prove that they are indeed rational in the following way:  Expand the series in powers of $x$ and compare with \cite[Expansion 1.1]{Gui-matrix} using \cite[eq. 2.30 \& 2.32]{Gui-matrix}. Below, we prove (\ref{bilateral-rama}).
\begin{proof}
The function $f: \mathbb{C} \longrightarrow \mathbb{C}$, defined by
\[
f(x)=e^{-i \pi x} (\cos \pi x)^{2j+1} \prod_{s_k \neq \frac12} \frac{\cos \pi x - \cos \pi s_k}{1-\cos \pi s_k}  \sum_{n=-\infty}^{\infty} A(n, x),
\]
is periodic of period $x=1$ because $e^{-i \pi x} (\cos \pi x)^{2j+1}$ is clearly periodic, the above product over $s_k \neq \frac12$ is periodic because each $s_k=s$, has a companion $s_k=1-s$, and the last factor is a sum of $A(n, x)=A(n+x, 0)$ over $\mathbb{Z}$, and therefore periodic as well. In addition $f(x)$ is holomorphic because the zero of $\cos \pi x - \cos \pi s_k$ at $x=-s_k$ cancels the pole of $(s_k)_{n+x}$ at $x=-s_k$, and as $f(x)$ is periodic all the other poles are canceled as well. As $f(x)$ is holomorphic and periodic, it has a Fourier expansion. Finally, we can prove that $f(x)=\mathcal{O}(e^{(2m+1) \pi |{\rm Im}(x)|})$, and therefore the Fourier expansion terminates at $k=m$.
\end{proof}
As (\ref{bilateral-rama}) must hold for all values of $x$, we can determine the values of $\alpha_k$ and $\beta_k$ by solving numerically the linear system of equations that results giving $2m$ values to $x$. Once we get the approximated solutions, is easy to identify the exact rational values they are.

\section{Bilateral semi-terminating series of Ramanujan-type}
If we let $x$ tend to $-p/2$, where $p$ is a positive odd number, and take into account that
\[
\lim_{x \to -\frac{p}{2}} \left(\frac12\right)_x (\cos \pi x) =  \frac{\sqrt{\pi}}{(\frac{p-1}{2})!} = \sqrt{\pi} \frac{(-1)^{\frac{p-1}{2}}} {\left( \frac12-\frac{p}{2}\right)_{\frac{p-1}{2}}}, \qquad (s_k)_{-\frac{p}{2}} = \frac{(s_k)_{-\frac12}}{\left(s_k-\frac{p}{2}\right)_{\frac{p-1}{2}}} \, \, \, \text{if} \, \, \, s_k \neq \frac12,
\]
we deduce the following bilateral identity:
\begin{multline}
\sum_{n=-\infty}^{\infty}  B\left(n, -\frac{p}{2}\right) =\sum_{n=-\infty}^{\frac{p-1}{2}} B\left(n, -\frac{p}{2}\right) 
= v_0 \frac{1+\sum_{k=1}^{m} \alpha_k ((-1)^k-1)}{\sum_{k=0}^{m} a_k \left( -\frac12\right)^k}
\\ \times \prod_{s_k \neq \frac12} \left( \frac{\Gamma(s_k)}{\Gamma(s_k - \frac12)} \, \frac{1-\cos \pi s_k}{\cos \pi s_k} \right) \cdot \frac{ \sqrt{ (-z_0) (-1)^{m} \chi_0}}{\pi^j} B\left(\frac{p-1}{2}, -\frac{p}{2}\right) \\ =
r_0 \, \frac{\sqrt{(-1)^{j} \varepsilon_0}}{\pi^j} \, B\left(\frac{p-1}{2}, -\frac{p}{2}\right),
\end{multline}
where  $\varepsilon_0$ is the discriminant of a quadratic field (imaginary or real), and we conjecture that $r_0$ is rational (because we have previously conjectured that all $\alpha_k$ and $\beta_k$ are rational). In addition, we conjecture the following supercongruences: 
\begin{equation}
\sum_{n=0}^{\frac{p-1}{2}} B\left(n, -\frac{p}{2}\right) \equiv \left( \frac{\varepsilon_0}{p} \right) p^j \, B\left(\frac{p-1}{2}, -\frac{p}{2}\right) \quad \pmod{p^{2m+1}},
\end{equation}
for odd primes $p$ (except perhaps a finite number of them). See the conjectured values of $\varepsilon_0$ (only a few have been proved; see Section \ref{proved-eps}), in the tables at the end of the paper.
\par In a similar way, observe that if $x=p$, where $p$ is a positive integer, the sum 
\begin{equation}
\sum_{n=-\infty}^{\infty} \left( \prod_{i=0}^{2m} \frac{(s_i)_{n+p}}{(1)_{n+p}} \right) \sum_{k=0}^m a_k (n+p)^k z_0^{n+p} \, = \frac{v_0\sqrt{(-1)^m \chi_0}}{\pi^m}
\end{equation}
terminates to the left at $n=-p$, and we get
\begin{equation}
\sum_{n=-p}^{\infty} A(n,p) = \frac{v_0\sqrt{(-1)^m \chi_0}}{\pi^m},
\end{equation}
and inspired by the pattern of supercongruences at $x=-p/2$, we write
\[
\sum_{n=-p}^{\infty} A(n,p) = \frac{v_0}{a_0} \frac{\sqrt{(-1)^m \chi_0}}{\pi^m} A(-p, p), \quad A(-p, p)=a_0.
\]
Then, we again observe the same pattern of congruence
\begin{equation}\label{eq-cong-zud}
\sum_{n=-p}^{-1} A(n,p) \equiv A(-p, p) \, p^m \left(\frac{\chi_0}{p}\right) \quad \pmod{p^{2m+1}}, \qquad A(-p, p) = a_0,
\end{equation}
for primes $p$ (except a finite set). As
\begin{equation}
\sum_{n=-p}^{-1} A(n, p) = \sum_{n=-p}^{-1} A(n+p, 0) = \sum_{n=0}^{p-1} A(n,0),
\end{equation}
the supercongruences (\ref{eq-cong-zud}) are equivalent to the Zudilin-type supercongruences:
\begin{equation}
\sum_{n=0}^{p-1} A(n,0) \equiv  a_0 p^m \left(\frac{\chi_0}{p}\right) \quad \pmod{p^{2m+1}}.
\end{equation}
That is
\begin{equation}
\sum_{n=0}^{p-1} \left( \prod_{i=0}^{2m} \frac{(s_i)_n}{(1)_n} \right) \sum_{k=0}^m a_k n^k z_0^n \,  \equiv  a_0 \, p^m \left(\frac{\chi_0}{p}\right) \quad \pmod{p^{2m+1}},
\end{equation}
for primes $p$ except a finite set. Y. Zhao adds an extra term to get supercongruences module $p^{2m+2}$ \cite{Zhao-2}.

\section{Some unproved examples}
In the examples below, we define $A(n, x)$ and $B(n, x)$ as indicated in (\ref{Anx}) and (\ref{Bnx}). All the content of this section is conjectured.
\begin{example} \rm
Let
\[
\sum_{n=0}^{\infty} \frac{\left(\frac12\right)_n\left(\frac13\right)_n\left(\frac23\right)_n \left(\frac16\right)_n\left(\frac56\right)_n}{(1)_n^5} \left( 1930n^2+549n+45 \right) \left( - \frac{3^6}{4^6} \right)^n=\frac{384}{\pi^2}.
\]
In this case we have 
\[
m=2, \quad v_0=384, \quad \chi_0=1, \quad j=0, 
\]
and we get
\[
\alpha_1=-\frac{14}{3}, \quad \alpha_2=2, \quad \beta_1=\beta_2=0,
\]
If we take $x=-p/2$ in $B(n, x)$, where $p$ is an odd positive integer, then we get
\[
\sum_{n=-\infty}^{\frac{p-1}{2}} B\left(n,-\frac{p}{2}\right) = \frac{93}{253} \,  B\left( \frac{p-1}{2}, -\frac{p}{2}  \right),
\]
and 
\[
\sum_{n=0}^{\frac{p-1}{2}} B\left(n,-\frac{p}{2}\right) \equiv B\left( \frac{p-1}{2}, -\frac{p}{2}  \right) \quad \pmod{p^5},
\]
for primes $p>3$. 
\par On the other hand, if we take $x=p$ in $A(n, x)$, where $p$ is a positive integer, we get
\[
\sum_{n=-p}^{\infty} A(n,p) = \frac{384}{\pi^2},
\]
and 
\[
\sum_{n=0}^{p-1} A(n,0) = \sum_{n=-p}^{-1} A(n,p) \equiv p^2 A(-p, p) \quad \pmod{p^5}, \qquad A(-p, p) = 45,
\]
for odd primes $p$. 
\end{example}

\begin{example} \rm 
Let
\[
\sum_{n=0}^{\infty} \frac{\left(\frac12\right)_n \left( \frac18 \right)_n \left( \frac38\right)_n \left( \frac58 \right)_n \left( \frac78 \right)_n}{(1)_n^5}  (1920n^2+304n+15) \left( \frac{1}{7^4} \right)^n = \frac{28\sqrt{28}}{\pi^2}.
\]
In this case we have
\[
m=2, \quad v_0=28, \quad \chi_0=28, \quad j=0, 
\]
and we get
\[
\alpha_1=-\frac{79}{2}, \quad \alpha_2=\frac{23}{2}, \quad \beta_1=\frac52 \, i, \quad \beta_2=-\frac32 \, i.
\]
If we take $x=-p/2$ in $B(n, x)$, we get
\[
\sum_{n=-\infty}^{\frac{p-1}{2}} B\left(n,-\frac{p}{2}\right) = - \frac{30}{7^4} \sqrt{-7} \,  B\left( \frac{p-1}{2}, -\frac{p}{2}  \right),
\]
and 
\[
\sum_{n=0}^{\frac{p-1}{2}} B\left(n,-\frac{p}{2}\right) \equiv \left(\frac{-7}{p}\right) B\left( \frac{p-1}{2}, -\frac{p}{2}  \right) \quad \pmod{p^5},
\]
for primes $p>3$. 
\par On the other hand, if we take $x=p$ in $A(n, x)$, where $p$ is a positive integer, we get
\[
\sum_{n=-p}^{\infty} A(n,p) = \frac{28 \sqrt{28}}{\pi^2},
\]
and 
\[
\sum_{n=0}^{p-1} A(n,0) = \sum_{n=-p}^{-1} A(n,p) \equiv \left( \frac{28}{p} \right) p^2 A(-p, p) \quad \pmod{p^5}, \qquad A(-p, p) = 15.
\]
for odd primes $p \neq 7$. 
\end{example}

\begin{example} \rm
Let
\begin{equation*}
\sum_{n=0}^{\infty}\frac{\left(\frac12\right)_n^3\left(\frac13\right)_n\left(\frac23\right)_n}{(1)_n^5} \left( 28n^2+18n+3 \right) (-27)^n=\frac{6}{\pi^2}.
\end{equation*}
In this case we have 
\[
m=2, \quad v_0=6, \quad \chi_0=1, \quad j=1, 
\]
and we get
\[
\alpha_1=\frac13, \quad \alpha_2=\frac14, \quad \beta_1=\beta_2=0.
\]
If we take $x=-p/2$ in $B(n, x)$, where $p$ is an odd positive integer, then we get
\[
\sum_{n=-\infty}^{\frac{p-1}{2}} B\left(n,-\frac{p}{2}\right) = \frac{3}{\pi} \,  B\left( \frac{p-1}{2}, -\frac{p}{2}  \right),
\]
and 
\[
\sum_{n=0}^{\frac{p-1}{2}} B\left(n,-\frac{p}{2}\right) \equiv \left(\frac{-4}{p}\right) p \, B\left( \frac{p-1}{2}, -\frac{p}{2}  \right) \quad \pmod{p^5},
\]
for primes $p>3$. 
\par On the other hand, if we take $x=p$ in $A(n, x)$, where $p$ is a positive integer, we get
\[
\sum_{n=-p}^{\infty} A(n,p) = \frac{6}{\pi^2},
\]
and 
\[
\sum_{n=0}^{p-1} A(n,0) = \sum_{n=-p}^{-1} A(n,p) \equiv p^2 A(-p, p) \quad \pmod{p^5}, \qquad A(-p, p) = 3.
\]
for odd primes $p$. 
\end{example}

\begin{example} \rm
Let
\begin{equation*}
\sum_{n=0}^{\infty}\frac{\left(\frac12\right)_n^5\left(\frac13\right)_n\left(\frac23\right)_n \left(\frac14\right)_n\left(\frac34\right)_n}{(1)_n^9} \left( 4528n^4 + 3180n^3+972n^2+147n+9 \right) \left(-\frac{27}{256}\right)^n=\frac{768}{\pi^4},
\end{equation*}
discovered by Y. Zhao. In this case we have 
\[
m=4, \quad v_0=768, \quad \chi_0=1, \quad j=2,
\]
and we get
\[
\alpha_1=-\frac{25}{6}, \quad \alpha_2=\frac83, \quad \alpha_3=-\frac23, \quad \alpha_4=\frac16, \quad  \beta_1=\beta_2=\beta_3=\beta_4=0.
\]
If we take $x=-p/2$ in $B(n, x)$, where $p$ is an odd positive integer, then we get
\[
\sum_{n=-\infty}^{\frac{p-1}{2}} B\left(n,-\frac{p}{2}\right) = \frac{3}{\pi^2} \, B\left( \frac{p-1}{2}, -\frac{p}{2}  \right),
\]
and 
\[
\sum_{n=0}^{\frac{p-1}{2}} B\left(n,-\frac{p}{2}\right) \equiv p^2 \, B\left( \frac{p-1}{2}, -\frac{p}{2}  \right) \quad \pmod{p^9},
\]
for odd primes $p \neq 5$. 
\par On the other hand, if we take $x=p$ in $A(n, x)$, where $p$ is a positive integer, we get
\[
\sum_{n=-p}^{\infty} A(n,p) = \frac{768}{\pi^4},
\]
and 
\[
\sum_{n=0}^{p-1} A(n,0) = \sum_{n=-p}^{-1} A(n,p) \equiv p^4 A(-p, p) \quad \pmod{p^9}, \qquad A(-p, p) = 9.
\]
for odd primes $p$. 
\end{example}

\section{Some proved examples}\label{proved-eps}
The first paper with some proofs of A-supercongruences (Zudilin-type):
\begin{equation}
\sum_{n=0}^{p-1} A(n, 0) \equiv \left(\frac{\chi_0}{p}\right) p^m A(-p, p) \quad \pmod{p^{2m+1}}, \qquad A(-p, p)=a_0,
\end{equation}
was \cite{Zudilin-supercong}. In it Wadim Zudilin managed to use some WZ (Wilf-Zeilberger) pairs for that purpose. Related papers are: \cite{Gui-Zud-diver-rama} in which we proved the supercongruences corresponding to some ``divergent" examples, and  \cite[Apendix]{Gui-Zud-q-rama}, with proofs that consist in applying the method based in $q$-analogues explained in \cite{Guo-Zud}. Other related papers are \cite{Ling-Long} and \cite{ChDeSw}.
\\ Here, we will use hypergeometric identities obtained by the WZ method in \cite{Gui-more-hyper-iden} for proving a few B-supercongruences, that is supercongruences of the form 
\begin{equation}
\sum_{n=0}^{\frac{p-1}{2}} B\left(n, \frac{p-1}{2}\right) \equiv \left(\frac{\varepsilon_0}{p}\right) p^j B\left( \frac{p-1}{2}, -\frac{p}{2}\right) \quad \pmod{p^{2m+1}},
\end{equation}
corresponding to Ramanujan-like series of ranks $3$ and $5$. In our proofs we will use the following identities valid for odd numbers $p$:

\begin{align}
\left(\frac12-\frac{p}{2}\right)_{\frac{p-1}{2}} &=(-1)^{\frac{p-1}{2}}\left(\frac{p-1}{2}\right)!, \label{id1} \\
\left( \frac14-\frac{p}{2} \right)_{\frac{p-1}{2}}\left( \frac34-\frac{p}{2} \right)_{\frac{p-1}{2}}&=\frac{8}{8^p}\frac{(2p-1)!}{(p-1)!}, \label{id2} \\
\left( \frac13-\frac{p}{2} \right)_{\frac{p-1}{2}}\left( \frac23-\frac{p}{2} \right)_{\frac{p-1}{2}}&=\frac{4}{3^{\frac{3p-3}{2}} 4^p}\frac{(\frac{p-1}{2})! (3p-2)!}{(p-1)!(\frac{3p-3}{2})!}, \label{id3}
\\
\left(1-\frac{p}{2}\right)_{\frac{p-1}{2}} &=(-1)^{\frac{p-1}{2}} \frac{2 (p-1)!}{2^p( \frac{p-1}{2})!}, \label{id4} 
\end{align}
and
\begin{align}
\left( \frac32+\frac{p}{2}\right)_{\frac{p-1}{2}} &=\frac{p!}{( \frac{p+1}{2})!}, \label{id5} \\
\left(\frac12\right)_{\frac{p-1}{2}}&=\frac{(p-1)!}{2^{p-1}(\frac{p-1}{2})!}, \label{id6} \\
\left(\frac12-p\right)_{\frac{p-1}{2}}&=\frac{(2p-1)!(\frac{p-1}{2})!}{2^{p-1}p!(p-1)!}. \label{id7}
\end{align}

\begin{example} \rm
Let
\[
\sum_{n=0}^{\infty} \frac{(\frac12)_n^3}{(1)_n^3}(6n+1) \left( \frac14 \right)^n = \frac{4}{\pi},
\]
then, for odd primes $p$, the following supercongruences hold:
\[
\sum_{n=0}^{\frac{p-1}{2}} B\left(n, -\frac{p}{2}\right) \equiv p B\left(\frac{p-1}{2}, - \frac{p}{2} \right) \quad \pmod{p^3},
\]
that is
\begin{equation}
\sum_{n=0}^{\frac{p-1}{2}} \frac{(\frac12-\frac{p}{2})_n^3}{(1-\frac{p}{2})_n^3} \left(1+6(n-\frac{p}{2}) \right) \left(\frac14\right)^n \equiv
- \frac{4p}{2^p} \frac{(\frac12-\frac{p}{2})_{\frac{p-1}{2}}^3}{(1-\frac{p}{2})_{\frac{p-1}{2}}^3} \quad \pmod{p^3}.
\end{equation}
\end{example}

\begin{proof}
From \cite[eq. 21]{Gui-more-hyper-iden}, we obtain the following terminating hypergeometric identity:
\begin{equation}\label{termina-1}
\sum_{n=0}^{\frac{p-1}{2}} \frac{(\frac12-\frac{p}{2})_n^3}{(1-\frac{p}{2})_n^3} \left(1+6(n-\frac{p}{2}) \right) \left(\frac{1}{4}\right)^n =
-\frac{8p^3}{(p+1)^2} \sum_{n=0}^{\frac{p-1}{2}} \frac{(\frac12-\frac{p}{2})_n^2}{(\frac32+\frac{p}{2})_n^2}(-1)^n,
\end{equation}
that we can write as
\begin{multline}
\sum_{n=0}^{\frac{p-1}{2}} \frac{(\frac12-\frac{p}{2})_n^3}{(1-\frac{p}{2})_n^3} \left(1+6(n-\frac{p}{2}) \right) \left(\frac{1}{4}\right)^n \\ =- \frac{8p^3}{(p+1)^2} \frac{(\frac12-\frac{p}{2})_{\frac{p-1}{2}}^2}{(1-\frac{p}{2})_{\frac{p-1}{2}}^2} (-1)^{\frac{p-1}{2}} - \frac{8p^3}{(p+1)^2} \sum_{n=0}^{\frac{p-3}{2}} \frac{(\frac12-\frac{p}{2})_n^2}{(\frac32+\frac{p}{2})_n^2}(-1)^n. 
\end{multline}
If $p$ is a prime number we easily deduce that
\begin{equation}
\sum_{n=0}^{\frac{p-1}{2}} \frac{(\frac12-\frac{p}{2})_n^3}{(1-\frac{p}{2})_n^3} \left(1+6(n-\frac{p}{2}) \right) \left(\frac{1}{4}\right)^n 
\equiv - \frac{8p^3}{(p+1)^2} \frac{(\frac12-\frac{p}{2})_{\frac{p-1}{2}}^2}{(\frac32+\frac{p}{2})_{\frac{p-1}{2}}^2} (-1)^{\frac{p-1}{2}} \pmod{p^3}.
\end{equation}
Hence, we have to prove that
\begin{equation}
- \frac{4p}{2^p} \frac{(\frac12-\frac{p}{2})_{\frac{p-1}{2}}^3}{(1-\frac{p}{2})_{\frac{p-1}{2}}^3} \equiv - \frac{8p^3}{(p+1)^2} \frac{(\frac12-\frac{p}{2})_{\frac{p-1}{2}}^2}{(\frac32+\frac{p}{2})_{\frac{p-1}{2}}^2}(-1)^{\frac{p-1}{2}} \quad \pmod{p^3}.
\end{equation}
Using the identities (\ref{id1}), (\ref{id4}), and (\ref{id5}), we see that it is a consequence of the Morley's congruence \cite{Morley}:
\begin{equation}\label{Morley}
\binom{p-1}{\frac{p-1}{2}} \equiv (-1)^{\frac{p-1}{2}} 2^{2(p-1)} \quad \pmod{p^3}, 
\end{equation}
and we are done.
\end{proof}
\begin{remark} \rm We can prove terminating identities like (\ref{termina-1}) in a direct an automatic way using the Zeilberger's algorithm for finding recurrences: Let $p=2k+1$ and denote with $\ell_k$ and $r_k$ the left and right hand sides respectively. Using the Zeilberger's algorithm we see that $\ell_k$ and $r_k$ satisfy the same recurrences. In addition, these recurrences are of first order. Finally, as both series are terminating, we can easily check that the necessary initial value is equal, that is $\ell_0 = r_0$. We can apply this method to all the examples in this section.
\end{remark}

\begin{example} \rm
Let
\[
\sum_{n=0}^{\infty} \frac{(\frac12)_n^3}{(1)_n^3}(42n+5) \left( \frac{1}{64} \right)^n = \frac{16}{\pi},
\]
then, for odd primes $p$, the following supercongruences hold:
\[
\sum_{n=0}^{\frac{p-1}{2}} B\left(n, -\frac{p}{2}\right) \equiv p B\left(\frac{p-1}{2}, - \frac{p}{2} \right) \quad \pmod{p^3},
\]
that is
\begin{equation}
\sum_{n=0}^{\frac{p-1}{2}} \frac{(\frac12-\frac{p}{2})_n^3}{(1-\frac{p}{2})_n^3} \left(5+42(n-\frac{p}{2}) \right) \left(\frac{1}{64}\right)^n \equiv - \frac{128p}{8^p} \frac{(\frac12-\frac{p}{2})_{\frac{p-1}{2}}^3}{(1-\frac{p}{2})_{\frac{p-1}{2}}^3} \quad \pmod{p^3}.
\end{equation}
\end{example}

\begin{proof}
From \cite[eq. 24]{Gui-more-hyper-iden}, we can obtain the following terminating hypergeometric identity:
\begin{multline}
\sum_{n=0}^{\frac{p-1}{2}} \frac{(\frac12-\frac{p}{2})_n^3}{(1-\frac{p}{2})_n^3} \left(5+42(n-\frac{p}{2}) \right) \left(\frac{1}{64}\right)^n =
-\frac{128p^3}{(p+1)^3} \sum_{n=0}^{\frac{p-1}{2}} \frac{(\frac12-\frac{p}{2})_n^3}{(\frac32+\frac{p}{2})_n^3}(-1)^n (2n+1) \\
=- \frac{128p^3}{(p+1)^3} \frac{(\frac12-\frac{p}{2})_{\frac{p-1}{2}}^2}{(\frac32+\frac{p}{2})_{\frac{p-1}{2}}^3} (-1)^{\frac{p-1}{2}} p - \frac{128p^3}{(p+1)^3} \sum_{n=0}^{\frac{p-3}{2}} \frac{(\frac12-\frac{p}{2})_n^3}{(\frac32+\frac{p}{2})_n^3}(-1)^n (2n+1).
\end{multline}
Hence
\begin{equation}
\sum_{n=0}^{\frac{p-1}{2}} \frac{(\frac12-\frac{p}{2})_n^3}{(1-\frac{p}{2})_n^3} \left(5+42(n-\frac{p}{2}) \right) \left(\frac{1}{64}\right)^n \equiv 
- \frac{128p^3}{(p+1)^3} \frac{(\frac12-\frac{p}{2})_{\frac{p-1}{2}}^3}{(\frac32+\frac{p}{2})_{\frac{p-1}{2}}^3} (-1)^{\frac{p-1}{2}} p \quad \pmod{p^3}.
\end{equation}
But, in view of the identities (\ref{id4}) and (\ref{id5}), we have
\begin{equation}
- \frac{128p^3}{(p+1)^3} \frac{(\frac12-\frac{p}{2})_{\frac{p-1}{2}}^3}{(\frac32+\frac{p}{2})_{\frac{p-1}{2}}^3} (-1)^{\frac{p-1}{2}} p = - \frac{128p}{8^p} \frac{(\frac12-\frac{p}{2})_{\frac{p-1}{2}}^3}{(1-\frac{p}{2})_{\frac{p-1}{2}}^3}.
\end{equation}
Hence, as they are equal, the congruences hold.
\end{proof}

\begin{example} \rm
Let
\[
\sum_{n=0}^{\infty} \frac{(\frac12)_n(\frac14)_n(\frac34)_n}{(1)_n^5}(20n+3) \left( -\frac14 \right)^n = \frac{8}{\pi},
\]
then, for odd primes $p$, the following supercongruences hold:
\[
\sum_{n=0}^{\frac{p-1}{2}} B\left(n, -\frac{p}{2}\right) \equiv  B\left(\frac{p-1}{2}, - \frac{p}{2} \right) \quad \pmod{p^3},
\]
that is
\begin{multline}
\sum_{n=0}^{\frac{p-1}{2}} \frac{(\frac12-\frac{p}{2})_n (\frac14-\frac{p}{2})_n (\frac34-\frac{p}{2})_n }{(1-\frac{p}{2})_n^3} \left(3+20(n-\frac{p}{2}) \right) \left(\frac{-1}{4}\right)^n \\ \equiv - \frac{14(-1)^{\frac{p-1}{2}}}{2^p} \frac{(\frac12-\frac{p}{2})_{\frac{p-1}{2}} (\frac14-\frac{p}{2})_{\frac{p-1}{2}} (\frac34-\frac{p}{2})_{\frac{p-1}{2}} }{(1-\frac{p}{2})_{\frac{p-1}{2}}^3} \quad \pmod{p^3}.
\end{multline}
\end{example}

\begin{proof}
Let $p$ be a positive odd number, replacing $x$ with $-p/2$ in \cite[eq. 27]{Gui-more-hyper-iden}, and taking into account that
\begin{equation}
8 \lim_{x \to \frac{-p}{2}} \left( \frac{4^x}{\pi} \, \frac{\cos 2\pi x}{\cos \pi x} \, \frac{(1)_x^3}{(\frac12)_x(\frac14)_x(\frac34)_x} \right)
= 2(-1)^{\frac{p-1}{2}}2^{-p} (2p-1)! \frac{(\frac{p-1}{2})!^4}{(p-1)!^4},
\end{equation}
we obtain the following terminating hypergeometric identity:
\begin{multline}
\sum_{n=0}^{\frac{p-1}{2}} \frac{(\frac12-\frac{p}{2})_n (\frac14-\frac{p}{2})_n (\frac34-\frac{p}{2})_n }{(1-\frac{p}{2})_n^3} \left(3+20(n-\frac{p}{2}) \right) \left(\frac{-1}{4}\right)^n = \\
2(-1)^{\frac{p-1}{2}}2^{-p} (2p-1)! \frac{(\frac{p-1}{2})!^4}{(p-1)!^4}
-\frac{32p^3}{p+1} \sum_{n=0}^{\frac{p-1}{2}} \frac{(\frac12-\frac{p}{2})_n^2 (\frac12-p)_n}{(\frac32+\frac{p}{2})_n(\frac12)_n^2} \frac{(-1)^n (2n+1-\frac{p}{2})}{(2n+1)^2} \\ = 
2(-1)^{\frac{p-1}{2}}2^{-p} (2p-1)! \frac{(\frac{p-1}{2})!^4}{(p-1)!^4} - 
(-1)^{\frac{p+1}{2}} \frac{16p^2}{p+1} \frac{(\frac12-\frac{p}{2})_{\frac{p-1}{2}}^2 (\frac12-p)_{\frac{p-1}{2}}}{(\frac32+\frac{p}{2})_{\frac{p-1}{2}}(\frac12)_{\frac{p-1}{2}}^2} \\
 -\frac{32p^3}{p+1} \sum_{n=0}^{\frac{p-3}{2}} \frac{(\frac12-\frac{p}{2})_n^2 (\frac12-p)_n}{(\frac32+\frac{p}{2})_n(\frac12)_n^2} \frac{(-1)^n (2n+1-\frac{p}{2})}{(2n+1)^2}
\end{multline}
If $p$ is a prime number we deduce that
\begin{multline}
\sum_{n=0}^{\frac{p-1}{2}} \frac{(\frac12-\frac{p}{2})_n (\frac14-\frac{p}{2})_n (\frac34-\frac{p}{2})_n }{(1-\frac{p}{2})_n^3} \left(3+20(n-\frac{p}{2}) \right) \left(\frac{-1}{4}\right)^n \equiv \\
2(-1)^{\frac{p-1}{2}}2^{-p} (2p-1)! \frac{(\frac{p-1}{2})!^4}{(p-1)!^4} - 
(-1)^{\frac{p-1}{2}} \frac{16p^2}{p+1} \frac{(\frac12-\frac{p}{2})_{\frac{p-1}{2}}^2 (\frac12-p)_{\frac{p-1}{2}}}{(\frac32+\frac{p}{2})_{\frac{p-1}{2}}(\frac12)_{\frac{p-1}{2}}^2} \pmod{p^3}.
\end{multline}
Then, using (\ref{id1}), (\ref{id2}), (\ref{id4}), (\ref{id5}), (\ref{id6}) and (\ref{id7}), we see that we have to prove that
\begin{multline}
2(-1)^{\frac{p-1}{2}}2^{-p}\frac{(2p-1)!}{p!(p-1)!}\frac{(\frac{p-1}{2})!^4}{(p-1)!^2} p - 4 \cdot 2^p \frac{(\frac{p-1}{2})^6}{(p-1)!^3} \frac{(2p-1)!}{p!(p-1)!}p 
\\ \equiv
-14(-1)^{\frac{p-1}{2}}2^{-p}\frac{(\frac{p-1}{2})!^4}{(p-1)!^2} \frac{(2p-1)!}{p!(p-1)!}p
\quad \pmod{p^3},
\end{multline}
or equivalently 
\begin{equation}
16(-1)^{\frac{p-1}{2}}p \, 2^{-p} \frac{\binom{2p-1}{p-1}}{\binom{p-1}{\frac{p-1}{2}}} 
\left( \binom{p-1}{\frac{p-1}{2}} -2^{2(p-1)} (-1)^{\frac{p-1}{2}} \right) \equiv 0 \quad \pmod{p^3},
\end{equation}
which is a consequence of the Morley's congruence (\ref{Morley}).
\end{proof}

\begin{example} \rm
Let 
\[
\sum_{n=0}^{\infty} \frac{(\frac12)_n^3(\frac14)_n(\frac34)_n}{(1)_n^5}(120n^2+34n+3) \left( \frac{1}{16} \right)^n = \frac{32}{\pi^2},
\]
then, for primes $p>3$, the following supercongruences hold:
\[
\sum_{n=0}^{\frac{p-1}{2}} B\left(n, -\frac{p}{2}\right) \equiv p B\left(\frac{p-1}{2}, - \frac{p}{2} \right) \quad \pmod{p^5},
\]
that is
\begin{multline}
\sum_{n=0}^{\frac{p-1}{2}} \frac{(\frac12-\frac{p}{2})_n^3 (\frac14-\frac{p}{2})_n(\frac34-\frac{p}{2})_n}{(1-\frac{p}{2})_n^5} \left(3+34(n-\frac{p}{2})+120(n-\frac{p}{2})^2 \right) \left(\frac{1}{2} \right)^{4n} \equiv \\ 
p \frac{64}{4^p}  \frac{(\frac12-\frac{p}{2})_{\frac{p-1}{2}}^3 (\frac14-\frac{p}{2})_{\frac{p-1}{2}}(\frac34-\frac{p}{2})_{\frac{p-1}{2}}}{(1-\frac{p}{2})_{\frac{p-1}{2}}^5} \quad \pmod{p^5}.
\end{multline}
\end{example}

\begin{proof}
From \cite[eq. 35]{Gui-more-hyper-iden}, we can obtain the following terminating hypergeometric identity:
\begin{multline}
\sum_{n=0}^{\frac{p-1}{2}} \frac{(\frac12-\frac{p}{2})_n^3 (\frac14-\frac{p}{2})_n(\frac34-\frac{p}{2})_n}{(1-\frac{p}{2})_n^5} \left(3+34(n-\frac{p}{2})+120(n-\frac{p}{2})^2 \right) \left(\frac{1}{2} \right)^{4n} = \\
\frac{256 p^5}{(p+1)^4} \sum_{n=0}^{\frac{p-1}{2}} \frac{(\frac12-\frac{p}{2})_n^4}{(\frac32+\frac{p}{2})_n^4}(2n+1) = 
\frac{256 p^5}{(p+1)^4} \frac{(\frac12-\frac{p}{2})_{\frac{p-1}{2}}^4}{(\frac32+\frac{p}{2})_{\frac{p-1}{2}}^4} p +\frac{256 p^5}{(p+1)^4} \sum_{n=0}^{\frac{p-3}{2}} \frac{(\frac12-\frac{p}{2})_n^4}{(\frac32+\frac{p}{2})_n^4}(2n+1).
\end{multline}
If $p$ is a prime number we easily deduce that
\begin{multline}
\sum_{n=0}^{\frac{p-1}{2}} \frac{(\frac12-\frac{p}{2})_n^3 (\frac14-\frac{p}{2})_n(\frac34-\frac{p}{2})_n}{(1-\frac{p}{2})_n^5} \left(3+34(n-\frac{p}{2})+120(n-\frac{p}{2})^2 \right) \left(\frac{1}{2} \right)^{4n} \equiv \\
\frac{256 p^5}{(p+1)^4} \frac{(\frac12-\frac{p}{2})_{\frac{p-1}{2}}^4}{(\frac32+\frac{p}{2})_{\frac{p-1}{2}}^4} p \quad \pmod{p^5}.
\end{multline}
We have to prove that
\begin{equation}
\frac{256 p^5}{(p+1)^4} \frac{(\frac12-\frac{p}{2})_{\frac{p-1}{2}}^4}{(\frac32+\frac{p}{2})_{\frac{p-1}{2}}^4} p  - \frac{64}{4^p}  \frac{(\frac12-\frac{p}{2})_{\frac{p-1}{2}}^3 (\frac14-\frac{p}{2})_{\frac{p-1}{2}}(\frac34-\frac{p}{2})_{\frac{p-1}{2}}}{(1-\frac{p}{2})_{\frac{p-1}{2}}^5}p \equiv 0 \quad \pmod{p^5}.
\end{equation}
Using the identities (\ref{id1}), (\ref{id2}), (\ref{id4}), and (\ref{id5}), we see that the above congruences are a consequence of 
\begin{equation}
\binom{2p-1}{p} \equiv 1 \quad \pmod{p^3},
\end{equation}
which is a result due to J. Wolstenholme \cite{Wolstenholme}.
\end{proof}

\begin{example} \rm
Let
\[
\sum_{n=0}^{\infty} \frac{(\frac12)_n^3(\frac13)_n(\frac23)_n}{(1)_n^5}(74n^2+27n+3) \left( \frac{27}{64} \right)^n = \frac{48}{\pi^2},
\]
then, for primes $p>3$, the following supercongruences hold:
\[
\sum_{n=0}^{\frac{p-1}{2}} B\left(n, -\frac{p}{2}\right) \equiv p B\left(\frac{p-1}{2}, - \frac{p}{2} \right) \quad \pmod{p^5},
\]
that is
\begin{multline}
\sum_{n=0}^{\frac{p-1}{2}} \frac{(\frac12-\frac{p}{2})_n^3 (\frac13-\frac{p}{2})_n(\frac23-\frac{p}{2})_n}{(1-\frac{p}{2})_n^5} \left(3+27(n-\frac{p}{2})+74(n-\frac{p}{2})^2 \right) \left(\frac{3}{4} \right)^{3n} \\ \equiv 8 p \left( \frac34 \right)^{3 \frac{p-1}{2}} \frac{(\frac12-\frac{p}{2})_{\frac{p-1}{2}}^3 (\frac13-\frac{p}{2})_{\frac{p-1}{2}}(\frac23-\frac{p}{2})_{\frac{p-1}{2}}}{(1-\frac{p}{2})_{\frac{p-1}{2}}^5} 
\end{multline}
\end{example}

\begin{proof}
From \cite[eq. 36]{Gui-more-hyper-iden}, we get the following terminating hypergeometric identity:
\begin{multline}
\sum_{n=0}^{\frac{p-1}{2}} \frac{(\frac12-\frac{p}{2})_n^3 (\frac13-\frac{p}{2})_n(\frac23-\frac{p}{2})_n}{(1-\frac{p}{2})_n^5} \left(3+27(n-\frac{p}{2})+74(n-\frac{p}{2})^2 \right) \left(\frac{3}{4} \right)^{3n} = \\
\frac{64 p^5}{(p+1)^3} \sum_{n=0}^{\frac{p-1}{2}} \frac{(\frac12-\frac{p}{2})_n^3}{(\frac32+\frac{p}{2})_n^3}=
\frac{64 p^5}{(p+1)^3} \frac{(\frac12-\frac{p}{2})_{\frac{p-1}{2}}^3}{(\frac32+\frac{p}{2})_{\frac{p-1}{2}}^3}+
\frac{64 p^5}{(p+1)^3} \sum_{n=0}^{\frac{p-3}{2}} \frac{(\frac12-\frac{p}{2})_n^3}{(\frac32+\frac{p}{2})_n^3}.
\end{multline}
If $p$ is a prime number we easily deduce that
\begin{multline}
\sum_{n=0}^{\frac{p-1}{2}} \frac{(\frac12-\frac{p}{2})_n^3 (\frac13-\frac{p}{2})_n(\frac23-\frac{p}{2})_n}{(1-\frac{p}{2})_n^5} \left(3+27(n-\frac{p}{2})+74(n-\frac{p}{2})^2 \right) \left(\frac{3}{4} \right)^{3n} \equiv \\
\frac{64 p^5}{(p+1)^3} \frac{(\frac12-\frac{p}{2})_{\frac{p-1}{2}}^3}{(\frac32+\frac{p}{2})_{\frac{p-1}{2}}^3} \quad \pmod{p^5}.
\end{multline}
We have to prove that
\begin{equation}
8 p \left( \frac34 \right)^{3 \frac{p-1}{2}} \frac{(\frac12-\frac{p}{2})_{\frac{p-1}{2}}^3 (\frac13-\frac{p}{2})_{\frac{p-1}{2}}(\frac23-\frac{p}{2})_{\frac{p-1}{2}}}{(1-\frac{p}{2})_{\frac{p-1}{2}}^5} \equiv \frac{64 p^5}{(p+1)^3} \frac{(\frac12-\frac{p}{2})_{\frac{p-1}{2}}^3}{(\frac32+\frac{p}{2})_{\frac{p-1}{2}}^3} \quad \pmod{p^5},
\end{equation}
Using the identities (\ref{id1}), (\ref{id3}), (\ref{id4}), and (\ref{id5}), we see that the above congruences are a consequence of the following ones:
\begin{equation}
p \dbinom{\frac{3p-3}{2}}{\frac{p-1}{2}, \frac{p-1}{2}, \frac{p-1}{2}} \equiv (-1)^{\frac{p-1}{2}}(3p-2)\dbinom{3p-3}{p-1, p-1, p-1}  \quad \pmod{p^5},
\end{equation}
for primes $p \geq 5$, which comes from \cite[eq. 19]{Zudilin-supercong} and Morley's congruence (\ref{Morley}).
\end{proof}

\begin{example} \rm
Let (see Remark \ref{obser-1}):
\[
\sum_{n=0}^{\infty} \frac{(\frac12)_n^5}{(1)_n^5}(10n^2+6n+1) (-4)^n = \frac{4}{\pi^2},
\]
then, for primes $p>3$, the following supercongruences hold:
\[
\sum_{n=0}^{\frac{p-1}{2}} B\left(n, -\frac{p}{2}\right) \equiv p^2 B\left(\frac{p-1}{2}, - \frac{p}{2} \right) \quad \pmod{p^5},
\]
that is
\begin{equation}
\sum_{n=0}^{\frac{p-1}{2}} \frac{(\frac12-\frac{p}{2})_n^5}{(1-\frac{p}{2})_n^5} \left(1+6(n-\frac{p}{2})+10(n-\frac{p}{2})^2 \right) (-4)^n \\ \equiv \frac14 \, p^2 \, (-1)^{\frac{p-1}{2}} \, 2^p \frac{(\frac12-\frac{p}{2})_{\frac{p-1}{2}}^5}{(1-\frac{p}{2})_{\frac{p-1}{2}}^5} 
\end{equation}
\end{example}
\begin{proof}
If $p$ is a positive odd number, then the following terminating hypergeometric identy holds
\begin{multline}
\sum_{n=0}^{\frac{p-1}{2}} \frac{(\frac12-\frac{p}{2})_n^5}{(1-\frac{p}{2})_n^5} \left(1+6(n-\frac{p}{2})+10(n-\frac{p}{2})^2 \right) (-4)^n =
p^5 \sum_{n=0}^{\frac{p-1}{2}} \frac{(\frac12-\frac{p}{2})_n^4}{(\frac12)_n^4} \frac{2n+1-\frac{p}{2}}{(2n+1)^4} 
\\ = \frac{p^6}{2p+2} \frac{(\frac12-\frac{p}{2})_{\frac{p-1}{2}}^4}{(\frac12)_{\frac{p-1}{2}}^4} + p^5 \sum_{n=0}^{\frac{p-3}{2}} \frac{(\frac12-\frac{p}{2})_n^4}{(\frac12)_n^4} \frac{2n+1-\frac{p}{2}}{(2n+1)^4}.
\end{multline}
If $p$ is a prime number we easily deduce that
\[
\sum_{n=0}^{\frac{p-1}{2}} \frac{(\frac12-\frac{p}{2})_n^5}{(1-\frac{p}{2})_n^5} \left(1+6(n-\frac{p}{2})+10(n-\frac{p}{2})^2 \right) (-4)^n \equiv
\frac{p^6}{2p+2} \frac{(\frac12-\frac{p}{2})_{\frac{p-1}{2}}^4}{(\frac12)_{\frac{p-1}{2}}^4} \quad \pmod{p^5}.
\]
Hence, we have to prove that
\[
\frac{p^2}{4} \, (-1)^{\frac{p-1}{2}} \, 2^p \frac{(\frac12-\frac{p}{2})_{\frac{p-1}{2}}^5}{(1-\frac{p}{2})_{\frac{p-1}{2}}^5} \equiv 
\frac{p^2}{2} \frac{(\frac12-\frac{p}{2})_{\frac{p-1}{2}}^4}{(\frac12)_{\frac{p-1}{2}}^4} \quad \pmod{p^5}.
\]
Using (\ref{id1}), (\ref{id4}) and (\ref{id6}), we see that it is a consequence of Morley's congruence (\ref{Morley}).
\end{proof}

\section{Tables}
In this section, we give tables of rational Ramanujan-like series for $1/\pi^m$ of ranks $2m+1=3, 5, 7, 9$. In them, we show the values of $z_0$, $a_0, a_1,...,a_m$, $v_0$, and the discriminant $\chi_0$ and $\varepsilon_0$. Tables \ref{3F2-1} and \ref{3F2-2} are complete, and in the other tables we show all the known series (most of them remain unproved).

\begin{table}[ht]
    \begin{tabular}{|c || c | c | c  c  c | c  c |}
        \hline &&&&&&& \\[-0.5em]
        $\#$ & $\bf{s}$ & $z_0$ & $a_0$ & $a_1$ & $v_0$ & $\chi_0$ & $\varepsilon_0$  \\[0.5em]
        \hline \hline &&&&&&& \\
        $1$ & $\frac12, \frac12, \frac12$ & $-1$ & $1$ & $4$  & $1$ & $-4$ & $-4$ \\[0.5em]
        $2$ & $\frac12, \frac12, \frac12$ & $-1/2^3$ & $1$ & $6$ & $1$ &  $-8$ & $-4$ \\[0.5em]
        $3$ & $\frac12, \frac12, \frac12$ & $1/2^2$ & $1$ & $6$ & $2$ & $-4$ & $1$ \\[0.5em]
        $4$ & $\frac12, \frac12, \frac12$ & $1/2^6$ & $5$ & $42$ & $8$ & $-4$ & $1$ \\[0.5em]
        $5$ & $\frac12, \frac12, \frac12$ & $-2^3$  & $2$ & $6$ & $1$ & $-4$ & $-8$  \\[0.5em]
        $6$ & $\frac12, \frac12, \frac12$ & $2^2$ & $1$ & $3$ & $-2$ & $1$ & $-4$ \\[0.5em]
        $7$ & $\frac12, \frac12, \frac12$ & $2^6$ & $8$ & $21$ & $-4$ & $1$ & $-4$ \\[0.35em]
        \hline &&&&&&& \\[-0.65em]
        $8$ & $\frac12, \frac14, \frac34$ & $-1/4$ & $3$ & $20$ & $4$ & $-4$ & $1$ \\[0.5em]
        $9$ & $\frac12, \frac14, \frac34$ & $-(16/63)^2$ & $8$ & $65$ & $9$ &  $-7$ & $28$ \\[0.5em]
        $10$ & $\frac12, \frac14, \frac34$ & $-1/48$ & $9$ & $84$ & $16$ & $-3$ & $1$ \\[0.5em]
        $11$ & $\frac12, \frac14, \frac34$ & $-1/18^2$ & $23$ & $260$ & $36$ & $-4$ & $1$ \\[0.5em]
        $12$ & $\frac12, \frac14, \frac34$ & $-1/(5\cdot 72^2)$ & $205$ & $3220$ & $144$ & $-20$ & $1$ \\[0.5em]
        $13$ & $\frac12, \frac14, \frac34$ & $-1/882^2$ & $1123$ & $21460$ & $1764$ & $-4$ & $1$ \\[0.5em]
        $14$ & $\frac12, \frac14, \frac34$ & $32/81$ & $4$ & $28$ & $9$ & $-4$ & $-8$ \\[0.5em]
        $15$ & $\frac12, \frac14, \frac34$ & $1/3^2$ & $1$ & $8$ & $2$ & $-3$ & $-3$ \\[0.5em]
        $16$ & $\frac12, \frac14, \frac34$ & $1/3^4$ & $8$ & $80$ & $9$ & $-8$ & $1$ \\[0.5em]
        $17$ & $\frac12, \frac14, \frac34$ & $1/7^4$ & $27$ & $360$ & $49$ & $-3$ & $-3$ \\[0.5em]
        $18$ & $\frac12, \frac14, \frac34$ & $1/99^2$ & $19$ & $280$ &  $18$ & $-11$ & $-11$ \\[0.5em]
        $19$ & $\frac12, \frac14, \frac34$ & $1/99^4$ & $8824$ & $211120$ & $9801$ & $-8$ & $-8$ \\[0.5em]
        $20$ & $\frac12, \frac14, \frac34$ & $-16/9$ & $1$ & $5$ &  $1$ & $-3$ & $12$ \\[0.5em]
        $21$ & $\frac12, \frac14, \frac34$ & $(4/3)^4$ & $8$ & $35$ & $-18$ & $1$ & $8$ \\[0.35em]
        \hline
    \end{tabular}
    \vskip 0.25cm
    \caption{Rational Ramanujan-type series for $1/\pi$ (First part)}\label{3F2-1}
\end{table}

\begin{table}[ht]
    \begin{tabular}{|c || c | c | c  c  c | c  c |}
        \hline &&&&&&& \\[-0.5em]
        $\#$ & $\bf{s}$ & $z_0$ & $a_0$ & $a_1$ & $v_0$ & $\chi_0$ & $\varepsilon_0$  \\[0.5em] 
        \hline \hline &&&&&&& \\
        $22$ & $\frac12, \frac13, \frac23$ & $-9/16$ & $3$ & $15$ & $4$ & $-3$ & $1$ \\[0.5em]
        $23$ & $\frac12, \frac13, \frac23$ & $-1/16$ & $7$ & $51$ & $12$ & $-3$ & $1$ \\[0.5em]
        $24$ & $\frac12, \frac13, \frac23$ & $-1/80$ & $5$ & $45$ &  $4$ & $-15$ & $1$ \\[0.5em]
        $25$ & $\frac12, \frac13, \frac23$ & $-1/2^{10}$ & $106$ & $1230$ & $192$ & $-3$ & $1$ \\[0.5em]
        $26$ & $\frac12, \frac13, \frac23$ & $-1/3024$ & $182$ & $2310$ & $216$ & $-7$ & $1$ \\[0.5em]
        $27$ & $\frac12, \frac13, \frac23$ & $-1/500^2$ & $827$ & $14151$  & $1500$ & $-3$ & $1$ \\[0.5em]
        $28$ & $\frac12, \frac13, \frac23$ & $1/2$ & $1$ & $6$ & $3$ & $-3$ & $-8$ \\[0.5em]
        $29$ & $\frac12, \frac13, \frac23$ & $2/27$ & $16$ & $120$ & $27$ & $-4$ & $-8$ \\[0.5em]
        $30$ & $\frac12, \frac13, \frac23$ & $4/5^3$ & $8$ & $66$ & $15$ & $-3$ & $-20$ \\[0.5em]
        $31$ & $\frac12, \frac13, \frac23$ & $-2^2$ & $4$ & $15$ & $3$ & $-3$ & $1$ \\[0.5em]
        $32$ & $\frac12, \frac13, \frac23$ & $27/2$ & $3$ & $10$ & $-3$ & $1$ & $8$ \\[0.5em]
        $33$ & $\frac12, \frac13, \frac23$ & $27/16$ & $3$ & $11$ & $-12$ & $1$ & $1$ \\[0.35em]
        \hline &&&&&&& \\[-0.65em]
        $34$ & $\frac12, \frac16, \frac56$ & $-4^3/5^3$ & $8$ & $63$ & $5$ & $-15$ & $1$ \\[0.5em]
        $35$ & $\frac12, \frac16, \frac56$ & $-3^3/8^3$ & $15$ & $154$ & $16$ & $-8$ & $1$ \\[0.5em]
        $36$ & $\frac12, \frac16, \frac56$ & $-1/8^3$ & $25$ & $342$ & $16$ & $-24$ & $1$ \\[0.5em]
        $37$ & $\frac12, \frac16, \frac56$ & $-9/40^3$ & $279$ & $4554$ & $80$ & $-120$ & $1$ \\[0.5em]
        $38$ & $\frac12, \frac16, \frac56$ & $-1/80^3$ & $789$ & $16254$ & $640$ & $-15$ & $1$ \\[0.5em]
        $39$ & $\frac12, \frac16, \frac56$ & $-1/440^3$ & $10177$ & $261702$ & $880$ & $-1320$ & $1$ \\[0.5em]
        $40$ & $\frac12, \frac16, \frac56$ & $-1/53360^3$ & $13591409$ & $545140134$ &  $213440$ & $-40020$ & $1$ \\[0.5em]
        $41$ & $\frac12, \frac16, \frac56$ & $(3/5)^3$ & $6$ & $56$ & $5$ & $-20$ & $-4$ \\[0.5em]
        $42$ & $\frac12, \frac16, \frac56$ & $4/5^3$ & $6$ & $66$ & $5$ & $-15$ & $-4$ \\[0.5em]
        $43$ & $\frac12, \frac16, \frac56$ & $(2/11)^3$ & $40$ & $504$ & $11$ & $-132$ & $-8$ \\[0.5em]
        $44$ & $\frac12, \frac16, \frac56$ & $(4/85)^3$ & $432$ & $7182$ & $85$ &  $-255$ & $-4$ \\[0.5em]
        \hline
    \end{tabular}
    \vskip 0.25cm
    \caption{Rational Ramanujan-type series for $1/\pi$ (Second Part)}\label{3F2-2}
\end{table}

\clearpage

\begin{table}[ht]
    \begin{tabular}{|c || c | c | c  c  c  c | c  c |}
        \hline &&&&&&&& \\[-0.5em]
        $\#$ & $\bf{s}$ & $z_0$ & $a_0$ & $a_1$ & $a_2$ & $v_0$ & $\chi_0$ & $\varepsilon_0$  \\[0.5em] 
        \hline \hline &&&&&&&& \\
        $1$ & $\frac12, \frac12, \frac12, \frac12, \frac12$ & $-1/2^2$ & $1$ & $8$  & $20$ & $8$ &  $1$ & $1$ \\[0.5em]
        $2$ & $\frac12, \frac12, \frac12, \frac12, \frac12 $ & $-2^2$ & $1$ & $6$ & $10$ & $4$ & $1$ & $1$ \\[0.5em]
        $3$ & $\frac12, \frac12, \frac12, \frac12, \frac12 $ & $-1/2^{10}$ & $13$ & $180$ & $820$ & $128$ & $1$ & $1$ \\[0.5em]
        $4$ & $\frac12, \frac12, \frac12, \frac12, \frac12 $ & $-2^{10}$ & $32$ & $160$ & $205$ & $16$ & $1$ & $1$ \\[0.5em]
        $5$ & $\frac12, \frac12, \frac12, \frac13, \frac23 $ & $(3/4)^3$ & $3$ & $27$ & $74$ & $48$ & $1$ & $1$ \\[0.5em]
        $6$ & $\frac12, \frac12, \frac12, \frac13, \frac23 $ & $-3^3$ & $3$ & $18$ & $28$ & $6$ & $1$ & $-4$ \\[0.5em]
        $7$ & $\frac12, \frac15, \frac25, \frac35, \frac45 $ & $-5^5/2^8$ & $30$ & $245$ & $483$ & $80$ & $1$ & $1$ \\[0.5em]
        $8$ & $\frac12, \frac12, \frac12, \frac14, \frac34 $ & $1/2^4$ & $3$ & $34$ & $120$ & $32$ & $1$ & $1$ \\[0.5em]
        $9$ & $\frac12, \frac13, \frac23, \frac14, \frac34 $ & $-1/48$ & $5$ & $63$ & $252$ & $48$ &  $1$ & $1$ \\[0.5em]
        $10$ & $\frac12, \frac13, \frac23, \frac14, \frac34 $ & $-3^3/2^4$ & $9$ & $75$ & $172$ & $48$ & $1$ & $1$ \\[0.5em]
        $11$ & $\frac12, \frac13, \frac23, \frac16, \frac56 $ & $-(3/4)^6$ & $45$ & $549$ & $1930$ & $384$ & $1$ & $1$ \\[0.5em]
        $12$ & $\frac12, \frac13, \frac23, \frac16, \frac56 $ & $(3/5)^6$ & $36$ & $504$ & $2128$ & $375$ &  $1$ & $-4$ \\[0.5em]
        $13$ & $\frac12, \frac13, \frac23, \frac16, \frac56 $ & $-1/80^3$ & $29$ & $693$ & $5418$ & $128$ &  $5$ & $1$ \\[0.5em]
        $14$ & $\frac12, \frac14, \frac34, \frac16, \frac56 $ & $-1/2^{10}$ & $15$ & $278$ &  $1640$ & $128$ & $12$ & $1$ \\[0.5em]
        $15$ & $\frac12, \frac18, \frac38, \frac58, \frac78 $ & $1/7^4$ & $15$ & $304$ & $1920$ & $28$ &  $28$ & $-7$ \\[0.5em]
        \hline
    \end{tabular}
    \vskip 0.25cm
    \caption{Rational Ramanujan-like series for $1/\pi^2$}\label{table-5F4}
\end{table}

\begin{table}[ht]
    \begin{tabular}{|c || c | c | c c c c c | c  c |}
        \hline &&&&&&&&& \\[-0.5em]
        $\#$ & $\bf{s}$ & $z_0$ & $a_0$ & $a_1$ & $a_2$ & $a_3$ & $v_0$ & $\chi_0$ & $\varepsilon_0$ \\[0.5em] 
        \hline \hline &&&&&&&&& \\
        $1$ & $\frac12, \frac12, \frac12, \frac12, \frac12, \frac12, \frac12$ & $1/2^6$ & $1$ & $14$  & $76$ & $168$ & $16$ &  $-4$ & $1$ \\[0.5em]
        $2$ & $\frac12, \frac12, \frac12, \frac12, \frac12, \frac12, \frac12$ & $2^6$ & $4$ & $32$ & $88$ & $84$ & $-24$ & $1$ & $-4$ \\[0.5em]
        $3$ & $\frac12, \frac12, \frac12, \frac12, \frac12, \frac13, \frac23$ & $27/4$ & $3$ & $27$ & $84$ & $92$ & $48$ & $1$ & $1$ \\[0.5em]
        \hline
    \end{tabular}
    \vskip 0.25cm
    \caption{Rational Ramanujan-like series for $1/\pi^3$}\label{table-7F6}
\end{table}

\clearpage

\begin{table}[ht]
    \begin{tabular}{|c || c | c | c c c c c c | c  c |}
        \hline &&&&&&&&&& \\[-0.5em]
        $\#$ & $\bf{s}$ & $z_0$ & $a_0$ & $a_1$ & $a_2$ & $a_3$ & $a_4$ & $v_0$ & $\chi_0$ & $\varepsilon_0$ \\[0.5em]
        \hline \hline &&&&&&&&&& \\
        $1$ & $\frac12, \frac12, \frac12, \frac12, \frac12, \frac12, \frac12, \frac14, \frac34$ & $1/2^{12}$ & $21$ & $466$  & $4340$ & $20632$ & $43680$ & $2048$ & $1$ & $1$ \\[0.5em]
        $2$ & $\frac12, \frac12, \frac12, \frac12, \frac12, \frac13, \frac23, \frac14, \frac34$ & $-3^3/2^8$ & $9$ & $147$ & $972$ & $3180$ & $4528$ & $768$ & $1$ & $1$ \\[0.5em]
        $3$ & $\frac12, \frac12, \frac12, \frac12, \frac12, \frac15, \frac25, \frac35, \frac45$ & $-5^5/2^{10}$ & $30$ & $425$ & $2275$ & $5600$ & $5532$ & $1280$ & $1$ & $1$ \\[0.5em]
        \hline
    \end{tabular}
    \vskip 0.25cm
    \caption{Rational Ramanujan-like series for $1/\pi^4$}\label{table-9F8}
\end{table}

{\bf Intriguing observation:} The values of $\varepsilon_0$ in Table \ref{table-5F4} coincide with those shown in \cite[Table 5.2]{dembele-zud} corresponding to some important Asai $L$-functions of Hilbert modular forms. This looks intriguing to me. What could be an explanation?

\section{Acknowledgements}
I am very grateful to Wadim Zudilin for valuable feedback and stimulating discussions.

\end{document}